\documentstyle[12pt,leqno]{article}
\newcommand{\be}{{\bf E}}
\newcommand{\bee}{{\bf e}}

\newcommand{\bs}{{\bf S}}
\newcommand{\bx}{{\bf X}}
\newcommand{\by}{{\bf Y}}

\newcommand{\bt}{{\bf T}}
\newcommand{\bq}{{\bf Q}}
\newcommand{\bg}{{\bf G}}
\newcommand{\bo}{{\bf 0}}

\newcommand{\nod}{\noindent}
\newcommand{\has}{\hspace{1cm} a.s.}

\newtheorem{theorem}{Theorem}[section]

\newtheorem{corollary}{Corollary}[section]
\newtheorem{lemma}{Lemma}[section]

\def\oo{\infty}
\def\noo{n\to\oo}

\def\le{\left(}
\def\ri{\right)}

\def\kx{{\bf x}}
\def\ky{{\bf y}}
\def\kz{{\bf z}}
\def\kt{{\bf t}}
\def\bo{{\bf 0}}
\def\bee{{\bf e}}

\def\bs{{\bf S}}

\def\sumon{\sum_{k=1}^{n}}

\def\zd{{\cal Z}_d}
\def\z2{{\cal Z}_2}

\def\begg{\begin{equation}}
\def\endd{\end{equation}}
\def\bege{\begin{eqnarray}}
\def\ende{\end{eqnarray}}

\def\pe{{\bf P}}
\def\ep{{\varepsilon}}
\def\al{{\alpha}}
\def\ga{{\gamma}}
\def\ld{{\lambda_d}}
\def\la{{\lambda}}
\def\gd{{\gamma_d}}
\def\xdn{{\xi^{(d)}(n)}}

\begin{document}
\centerline{\Large\bf Maximal Local Time of a d-dimensional }
\bigskip
\centerline{\Large\bf Simple  Random Walk on Subsets.}

\bigskip
\bigskip
\bigskip
\bigskip
\bigskip

\renewcommand{\thefootnote}{1}
\noindent {\bf Endre Cs\'{a}ki$\!$} \footnote{Research supported
by the Hungarian National Foundation for Scientif\/ic Research,
Grant No. T 037886 and T 043037.}\\
Alfr\'ed R\'enyi Institute of Mathematics, Hungarian Academy
of Sciences, Budapest, P.O.B. 127, H-1364, Hungary. E-mail
address: csaki@renyi.hu

\bigskip

\renewcommand{\thefootnote}{2}
\noindent {\bf Ant\'{o}nia F\"{o}ldes$\,$}\footnote{Research
supported by a PSC CUNY Grant, No. 65685-0034. }\\
City University of New York, 2800 Victory Blvd., Staten Island,
New York 10314, U.S.A. E-mail address: afoldes@gc.cuny.edu

\bigskip

\noindent {\bf P\'al R\'ev\'esz}$^1$ \\
Institut f\"ur Statistik und Wahrscheinlichkeitstheorie,
Technische Universit\"at Wien, Wiedner Hauptstrasse
8-10/107 A-1040 Vienna, Austria. E-mail address:
reveszp@renyi.hu

\bigskip
\bigskip
\bigskip

\noindent {\it Abstract}: Strong theorems are given  for the
maximal local time on balls and subspaces for the $d$-dimensional
simple symmetric random walk.

\bigskip

\noindent  AMS 2000 Subject Classification: Primary
60J15; Secondary 60F15, 60J55.

\bigskip

\noindent Keywords: maximal local time, simple  random walk
in $d$-dimension, strong theorems. \vspace{.1cm}

\vfill
\newpage
\renewcommand{\thesection}{\arabic{section}.}
\section{Introduction and main results}
\renewcommand{\thesection}{\arabic{section}}
\setcounter{equation}{0}
\setcounter{theorem}{0}
\setcounter{lemma}{0}

\noindent Consider a simple symmetric random walk $\{\bs_n\}_{n=1}^{\oo} $
starting at the origin $\bo$ on the $d$-dimensional integer lattice ${\cal
Z}_d$, i.e. $\bs_0=\bo$, $\bs_n=\sumon \bx_k$, $n=1,2,\dots$, where
$\bx_k,\, k=1,2,\dots$ are i.i.d. random variables with distribution
$$
\pe (\bx_1=\bee_i)=\pe (\bx_1=-\bee_i)=\frac{1}{2d},\qquad i=1,2,...,d
$$
and $\{\bee_1,\bee_2,...\bee_d\}$ is a system of orthogonal unit
vectors in ${\cal Z}_d.$
 Define the local time of the walk by \begg
\xi^{(d)}(\kx,n):=\#\{k: \,\,0< k \leq n,\,\,\, \bs _k=\kx \}.
\label{loc1}
\endd
where $\kx$ is any lattice point of ${\cal Z}_d.$ The maximal
local time of the walk is defined as
 \begg \xi^{(d)}(n):=\max_{\kx \in {\cal Z}_d}\xi^{(d)}(\kx,n) .
\label{loc2}
\endd
The properties of $\xi^{d}(n)$ were extensively studied in the
cases $d=1$, $d=2$ and $d\geq 3$ separately. For $d=1$
the interested reader should consult the monograph of P.
R\'ev\'esz \cite{R90}. In this paper we are interested in
investigating the  maximum local time for $d\geq 2$ in a
restricted sense, namely we want to investigate the maximum on
certain subsets of the state space. It is easy to see that these
maximums depend on both of the size and the shape of the selected
subset. We will only investigate two types of subsets: balls
centered at the origin and subspaces.

\bigskip
\nod
{\bf 1.1 Two dimension.}

\medskip
\nod
 In what follows we present the most important results on
 local time for $d=2$ which are relevant to our investigation.
\bigskip
\newline
\nod {\bf Theorem A} (Erd\H{o}s and Taylor \cite{ET60})

 $$
\lim_{n \to \infty} \pe(\xi^{(2)}(\bo,n)<x\log n)=1-e^{-\pi x}.
$$
\nod
{\it Let $f(x)$  resp. $g(x)$ be a decreasing resp. increasing
 function for which
$f(x)\log x \uparrow \oo$, $g(x)(\log x)^{-1} \downarrow 0.$
Then
$$
\pi^{-1}g(n) \log n \leq\xi^{(2)}(\bo,n)
$$
finitely often with probability one if and only if
$$
\int_1^{\oo}\frac{g(x)}{x\log
x}e^{-g(x)}\,dx<\oo,
$$

$$
f(n)\log n\geq \xi^{(2)}(\bo,n)
$$
finitely often with probability one if and only if
$$
\int_1^{\oo}\frac{f(x)}{x\log x}\,dx<\oo.
$$
}

\bigskip
\nod {\bf Theorem B} (Erd\H{o}s and Taylor \cite{ET60})

$$
\frac{1}{4\pi}\leq \liminf_{\noo} \frac{\xi^{(2)}(n)}{(\log n)^2}\leq
\limsup_{\noo}
\frac{\xi^{(2)}(n)}{(\log n)^2}\leq \frac{1}{\pi}\has
$$

They also conjectured that the upper bound in the above theorem is
the correct limit. This was confirmed recently by Dembo {\it et al.}
\cite{DPRZ01}. In fact they proved the following more general result:

\bigskip
\nod {\bf Theorem C} (\cite{DPRZ01}) {\it Let $\bs_n=\sumon   \bx_k$
be an aperiodic  random walk with i.i.d. increments $ \bx_k \in
{\cal Z}_2$  that satisfy $\be\bx_1=0$ and $\be|\bx_1|^m <\oo$
for all $m<\oo.$ Denote by $\Gamma=\be\bx \bx'$ the covariance
matrix of the increments, and write $\pi_{\Gamma}:= 2\pi(det
\Gamma)^{1/2}.$ Define the local time and maximum local time as in
the simple walk case.  Let $M(n,\alpha)$ denote the number  of
points in the set $\{\kx:\xi(\kx,n)\geq \alpha(\log n)^2\}$. Then

$$\lim_{\noo} \frac{\xi^{(2)}(n)}{(\log n)^2}=\pi^{-1}_{\Gamma} \has$$
 and for $\alpha \in(0,\pi^{-1}_{\Gamma}]$

$$\lim_{\noo} \frac{\log M(n,\alpha)}{\log n}=1-\alpha \pi_{\Gamma}
\has$$
\nod
Moreover any {\rm(}random{\rm)} sequence $\{\kx_n\}$ such
that $\xi^{(2)}(\kx_n,n)/\xi^{(2)}(n) \to 1$, must satisfy}
 $$\lim_{\noo} \frac{\log|\kx_n|}{\log n}=\frac{1}{2} \has$$
For the simple symmetric random walk $\pi^{-1}_{\Gamma}=\pi.$

\bigskip
 As the above results show, the local time of every fixed point is roughly
 around $\log n$ but the maximal local time is around $(\log
 n)^2.$ This phenomenon suggests that taking the maximum  local time on
appropriate subsets,
 one might get orders in between $\log n$  and $(\log n)^2.$  The
 following result of Auer is in fact telling us that such a set needs
 to be quite big.

\bigskip
\noindent {\bf Theorem D} (Auer \cite{A90}) {\it For any $\ep>0$ we
have
$$\lim_{\noo}\sup_{||\kx||\leq
r_n}\left|\frac{\xi^{(2)}(\kx,n)}{\xi^{(2)}(\bo,n)}-1\right|=0 \has,$$
where
$$r_n=\exp((\log n)^{1/2}(\log\log n)^{-1/2-\ep}) $$
and $||\kx||$ stands for the usual Euclidean norm.}

\bigskip
Now let $A$ be a subset of ${\cal Z}_2$ and define

\begg
\xi^{(2)}_A(n):=\max_{\kx \in A}\xi^{(2)}(\kx,n).
\endd

\noindent Let moreover $B(r)$ denote the set of lattice points in the
disc of radius $r$ centered at the origin, i.e.
\begg
B(r):=\{\kx\in {\cal Z}_2:\, ||\kx||\leq r\}.
\endd

 Then with $r_n$ as in Theorem D, for any subset   $A \subseteq B(r_n),\,$
 $\xi^{(2)}_A(n)/\xi^{(2)}(\bo,n)\to 1$ as $n\to\infty$. Hence for
$\xi^{(2)}_A(n)$ we have the same asymptotic results as for
$\xi^{(2)}(\bo,n)$.

Denote by $L=L(a_1,a_2)$ the lattice points $\kx=(x_1,x_2)$ on the
line $a_1x_1+a_2x_2=0$, where $a_1$ and $a_2$ are integers not
both of them zero. Now we formulate our results for the
two-dimensional case. Our first theorem is telling us that for any
line going through the origin  which  contains lattice points at
all, the maximal local time has the same order of magnitude as for
the whole plane.
\begin{theorem}
For any line $L=L(a_1,a_2)$ such that $a_1,a_2$ are integers, not both of
them zero, we have
\begg
\frac{1}{8\pi}\leq
\liminf_{\noo}\frac{\xi^{(2)}_L(n)}{(\log n)^2}\leq \limsup_{\noo}
\frac{\xi^{(2)}_L(n)}{(\log n)^2}\leq\frac{1}{2\pi}\has
\label{line1}
\endd
\end{theorem}

The next two theorems contain our results about discs centered at
the origin.
\begin{theorem}
Let $r_n=n^{\al},\quad 0<\al\leq 1/2.$ Then
\begg \frac{4
\al^2}{\pi}\leq \liminf_{\noo}\frac{\xi^{(2)}_{B(r_n)}(n)}{(\log n)^2}\leq
\limsup_{\noo} \frac{\xi^{(2)}_{B(r_n)}(n)}{(\log n)^2}\leq
\frac{2\alpha}{\pi}\has \label{circ1}
\endd
\end{theorem}

\begin{theorem}
Let $r_n=\exp((\log n)^{\beta})$. For any $\ep>0,\quad  1/2\leq
\beta<1$, and large enough $n$ we have
 \begg \frac{4(1-\ep)}{\pi}(\log n)^{2\beta}\leq \xi^{(2)}_{B(r_n)}(n)\leq
  (\log n)^{2\beta+\ep} \has \label{circ2}\endd
\end{theorem}

\begin{corollary}
If  $L$ and $B(r_n)$ are the sets in Theorems {\rm 1.1, 1.2} and {\rm 1.3}
respectively, then  for $r_n=n^{\al},\quad 0<\al\leq 1/2$ we have

\begg \frac{
\al^2}{2\pi}\leq \liminf_{\noo}\frac{\xi^{(2)}_{B(r_n)\cap L}(n)}{(\log
n)^2}\leq
\limsup_{\noo} \frac{\xi^{(2)}_{B(r_n)\cap L}(n)}{(\log n)^2}\leq
\frac{1}{\pi}\min\le \frac12,\, 2\alpha\ri \has \label{circ3}\endd
and for $r_n=\exp((\log n)^{\beta})$ with $\quad  1/2\leq \beta <1$ we
have for any $\ep>0$ and large $n$

\begg \frac{1-\ep}{2\pi}(\log n)^{2\beta}\leq \xi^{(2)}_{B(r_n)\cap
L}(n)\leq
  (\log n)^{2\beta+\ep} \has \label{circ4}\endd
\end{corollary}

\bigskip
\nod
{\bf 1.2 Three and higher dimension.}

\medskip
\noindent Just like in two dimension, for a subset $A\subseteq {\cal Z}_d$
we define
\begg
\xi^{(d)}_A(n):=\max_{\kx \in A}\xi^{(d)}(\kx,n).
\endd

To formulate the most important  known results on $\xi^{(d)}(n)$ of
(1.2), we need some more definition. Denote by $\gamma_d(n)$ the
probability that in the first $n-1$ steps the path does not
return to the origin. Then
$$1=\ga_d(1)\geq \ga_d(2)\geq ...\geq \ga_d(n)\geq...>0.$$
It was proved in Chung and Hunt \cite{CH49} that, for $d\geq 3$
\begg \ga_d(n) \to\ga_d>0, \label{gam1}
\endd
and
 \begg
\ga_d<\ga_d(n)<\ga_d+O(n^{1-d/2}) \label{gam2}
\endd
as $\noo$. So $\gamma_d$ is the probability that the $d$-dimensional
simple symmetric random walk never returns to its starting point.

Let $\xi^{(d)}(\bo,\infty)$ be
the total local time at $0$ of the infinite path in $\zd$. Then (see
\cite{ET60}) $\,\xi^{(d)}(\bo,\infty)$ has geometric distribution:

\begg
\pe(\xi^{(d)}(\bo,\infty)=k)=\ga_d(1-\ga_d)^k,\qquad k=0,1,2,...
\label{geo}
\endd
Erd\H{o}s and Taylor \cite{ET60} proved the following strong law for the
maximal local time:

\bigskip
\nod
{\bf Theorem F} (\cite{ET60}){\it For $d\ge 3$
\begg
\lim_{\noo}\frac{\xdn}{\log n}=\ld \has, \label{la}
\endd
where}
\begg
\ld=-\frac{1}{\log(1-\gd)}.
\endd

\nod
 {\bf Remark.} For the exact value of $\gamma_3$ see e.g. Spitzer
\cite{S76}, p. 103 which implies that $\lambda_3 <1$ and hence $\ld<1$ for
all $d\geq 3.$

Let $B(r)$ stand for the (discrete) ball centered at the
origin in the $d$-dimensional space and having radius $r,$ i.e.

\begg
B(r):=\{\kx\in {\cal Z}_d:\, ||\kx||\leq r\}.
\endd

Let furthermore $\kx=(x_1,x_2,\dots,x_d)$,
$$ S_{d-1}:=\{\kx\in \zd:\, a_1x_1+ a_2x_2+...+a_dx_d  =0   \}
$$
and
$$S_{d-2}:=\{\kx\in\zd:\, a_1x_1+ a_2x_2+...+a_dx_d  =0, \quad b_1x_1+
b_2x_2+...+b_dx_d
=0\}
$$
with integer coefficients $a_1, a_2, ... a_d,  b_1, b_2, ...b_d.$

\bigskip
\noindent
 For subspaces we will prove the
following two results.
\begin{theorem} Suppose that $a_1, a_2, ... a_d$ are integers, not all of
them zero, then
$$
\lim_{\noo}\frac{\xi_{S_{d-1}}^{(d)}(n)}{\log n}=\frac{\ld}{2}
\has
$$
\end{theorem}

The above theorem is telling us that the maximal local time in the
$d-1$ dimensional subspace has the same order of magnitude as in
the whole $d$-dimensional space. On the other hand, the next theorem
 shows that in the $d-2$ dimensional subspace the maximal local
time gets drastically smaller.
\begin{theorem} Suppose that $a_1, a_2, ... a_d$ are integers, not all of
them zero and $b_1, b_2, ...b_d$ are also integers not all of them zero.
Assume also that the vectors $(a_1,a_2,\dots, a_d)$ and
$(b_1,b_2,\dots,b_d)$ are not parallel. Then

$$
\lim_{\noo}\frac{\xi_{S_{d-2}}^{(d)}(n)}{\log\log n}=\ld \has
$$
\end{theorem}

\bigskip

For balls centered at the origin we will prove the following
result:
\begin{theorem} For any  sequence $r_n \uparrow \oo,$ such that
$\limsup_{n\to\infty}(\log r_n)/(\log n)\leq 1/2$,
we have
\begg
\lim_{n\to\infty}\frac{\xi_{B(r_n)}^{(d)}(n)}{\log r_n}=2\ld \has
\endd
\end{theorem}

The organization of the paper is as follows. In Section 2 we
will present some well-known facts and prove some preliminary results.
Sections 3 and 4 contain the proofs of the two dimensional and higher
dimensional results, respectively. In Section 5 some implications of the
above results and some open questions will be discussed.
Throughout the paper $c,\, c_1,\dots, C,\, C_1,\dots$ will
denote positive constants, the value of which is unimportant and may
vary from line to line.

\renewcommand{\thesection}{\arabic{section}.}
\section{Preliminary facts and  results}
\renewcommand{\thesection}{\arabic{section}}
\setcounter{equation}{0} \setcounter{theorem}{0}
\setcounter{lemma}{0}

First we present the $d$-dimensional law of the iterated logarithm and
rate of escape for simple walk.

\bigskip

\noindent
{\bf Fact 1.} (Dvoretzky and Erd\H os \cite{DE50} or \cite{R90}, pp.
193, 195) {\it For a simple symmetric random walk in} $\zd$
\begg
\limsup_{\noo} (2n\log\log n)^{-1/2}d^{1/2}||\bs_n||=1 \has
\label{lil}
\endd

{\it Moreover, in case $d\geq 3$, for any $0<\ep<1/2$ and large enough $n$
we have}
\begg
||\bs_n||>n^{1/2-\ep}\qquad\has
\label{esc}
\endd

\medskip
Consider now the case $d=2$. Introduce

\begin{eqnarray}
p(\kx):&=& \pe(\min\{n: \bs_n=\bo\}>\min\{n: \bs_n=\kx\})
\nonumber\\&=&\pe(\{\bs_n\}\,\, {\rm reaches\,\, }\kx \,\,{\rm before
\,\,returning \,\,to}\,\, \bo).
 \label{pr1}
\end{eqnarray}

We will need the following two lemmas from R\'ev\'esz \cite{R90}:

\bigskip
\nod
 {\bf Fact 2.} (\cite{R90}, p. 219) {\it For a simple symmetric random
walk in $\z2$ there exists a positive constant $C$ such that for  any }
$\kx \in \z2$ {\it with} $||\kx||\geq 2$
 \begg
p(\kx)\geq \frac{C}{\log ||\kx||}.
 \endd

Let us define
 \begg
\rho_0:=0, \qquad \rho_n:=\min\{ k: k>\rho_{n-1},\, \bs _k=\bo  \},
\quad n=1,2,\dots
\label{rho}
\endd

\bigskip
 \noindent {\bf Fact 3.} (\cite{R90}, pp. 219-220.)
{\it For a simple symmetric random walk in $\z2$ let
 \begg
 Y_i(\kx):=\xi^{(2)}(\kx,\rho_i)-\xi^{(2)}
 (\kx,\rho_{i-1}), \quad i=1,2,\dots
\endd
Then for fixed $\kx\in\z2,\,$ $\{Y_i(\kx)\}_{i=1}^\infty$ are i.i.d.
random variables with the following distribution:}

\begg \pe(Y_1(\kx)=0)=1-p(\kx)\endd

\begg \pe(Y_1(\kx)=k)=(1-p(\kx))^{k-1} p^2(\kx), \quad k=1,2,\dots
\label{dist1}
\endd

\medskip
Now we will prove our first
\begin{lemma} For a simple symmetric random walk in $\z2$ we have
for arbitrary $\kx \in \z2$, and any $u>0$
\begg \pe\le\sum_{k=1}^n Y_k(\kx)>u\,n \right) \leq e^{np(\kx)(1-u/2)}.
\endd
\end{lemma}

\medskip\noindent
{\bf Remark.} Fact 3 and Lemma 2.1 are true for more general random walk,
but we need them only for simple symmetric case.

\medskip
\noindent {\bf Proof:} From (\ref{dist1}) we easily get with
$q(\kx)=1-p(\kx)\,$ that for any $z>0,$ for which $q(\kx)e^z <1,$

\begg \be\le
e^{zY_1(\kx)}\ri=q(\kx)+\frac{p^2(\kx)e^z}{1-q(\kx)e^z}.
\label{exp1}
\endd
Putting $z=\log(2/(1+q(\kx)))$, we have
\begg
q(\kx)e^z=\frac{2q(\kx)}{1+q(\kx)}<1 \label{q1}
\endd
so (\ref{exp1}) holds. Thus

 \begg
\be\le e^{zY_1(\kx)}\ri=
q(\kx)+\frac{p^2(\kx)e^z}{1-q(\kx)e^z}=1+p(\kx).
\label{exp2}
\endd
 By exponential Markov inequality, and (\ref{exp2})  we have
 \begin{eqnarray}
\pe\le\sum_{k=1}^n Y_k(\kx)>u\,n \right)&\leq&\frac{\le\be
e^{zY_1(\kx)}\ri^n}{e^{zun}}=(1+p(\kx))^n\le 1-\frac{p(\kx)}2\ri^{un}
\nonumber\\
&\leq&\exp\le np(\kx)(1-u/2)\ri,
\end{eqnarray}
where the inequality $1+v\leq e^v$ was used.  $\Box$

\medskip
For our next lemma we need further notations and facts.
Our main source for these is Spitzer's book \cite{S76}. Here we consider
a two-dimensional symmetric aperiodic recurrent walk on $\z2$, more
general than a simple symmetric random walk. Recall that a random walk in
$\zd$ is aperiodic if the steps are not supported on a proper subgroup of
$\zd$. All what we are quoting however work in case of a one-dimensional
walk under the same constraints as well. We adopt the notations and
definitions listed below from \cite{S76}. Denote the $n$-step probability
transition function by
$$
P_n(\kx,\ky)=P_n(\ky,\kx)=P_n(\bo,\,\kx-\ky)=
\pe(\bs_n=\kx-\ky|\bs_0=\bo ),\quad n=0,1,\dots,\quad \kx,\ky\in\z2.
$$

For $a\geq 0$ integer we define
\begg
\bt_a:=\min\{j>0: \bs_{a+j}=\bo\}
 \label{sp1}
\endd
and will denote $\bt_0=:\bt.$ Let (see \cite{S76}, pp. 107, 160-161)

\begg
\bq_n(\kx,\ky):=\pe_{\kx}(\bs_n=\ky,\, \bt>n),\quad n=0,1,\dots,\quad
\kx,\ky\in\z2-\{\bo\},
 \label{sp2}
\endd

\noindent and
\begg
g(\kx,\ky):=\sum_{n=0}^{\oo} \bq_n(\kx,\ky), \quad \kx,\ky\in\z2-\{\bo\}
\label{sp3}
\endd
where $\pe_{\kx}(\cdot):=\pe(\cdot|\bs_0=\kx)$.
Using $\pe(\cdot):=\pe_{\bo}(\cdot)$ we also
recall from Spitzer \cite{S76} that
\begg
\pe_{\kx}(\bt>n)=\sum_{\kt\neq \bo}g(\kx, \kt)\sum_{\ky\neq
\bo}\bq_n(\kt, \ky)P(\ky,\bo), \quad n=0,1,\dots,\, \kx\in\z2-\{\bo\}
\label{sp4}
\endd
and introducing
\begg v_n(\kt)=\frac{1}{\pe(\bt>n+1)} \sum_{\ky\neq
\bo}\bq_n(\kt,\ky)P(\ky,\bo), \quad n=0,1,\dots,\,  \kt \in \z2-\{\bo\},
\label{sp5}
\endd

\noindent we get for $\kx\neq \bo$
\begg
\pe_{\kx}(\bt>n)=\pe(\bt>n+1)\sum_{\kt \neq \bo}g(\kx,\kt)v_n(\kt),
\quad n=0,1,\dots
 \label{sp6}
\endd
with \begg v_n(\kt)\geq0, \qquad \sum_{\kt \neq \bo}v_n(\kt)=1.
\label{sp7}
\endd

We introduce further the notations (see \cite{S76}, pp.
114, 139, 328)
\begg \bg_n(\kx,\ky)= \sum_{k=0}^n P_k(\kx, \ky),
\quad n=0,1,\dots,\, \kx,\ky\in\z2,
\label{sp8}
\endd
the truncated Green function
\begg g(n)=\bg_n(\bo,\bo)=
\sum_{k=0}^n P_k(\bo, \bo),\quad n=0,1,\dots \label{sp9}
\endd

\noindent and the potential kernel
\begg
a(\kx)=\sum_{n=0}^{\oo}(P_n(\bo,\bo)-P_n(\bo,\kx)),\quad \kx\in\z2.
\label{sp10}
\endd
Then it is known, that

\begg 0\leq g(\kx,\ky)=a(\kx)+a(\ky)-a(\kx-\ky). \label{sp11}
\endd

We recall from Spitzer \cite{S76}, p. 139

\bigskip
\noindent {\bf Fact 4.} {\it For any symmetric recurrent aperiodic
random walk in two dimension}
 \begg
\sum_{\ky\in\z2}P_{n+1}(\kx,\ky)a(\ky)=a(\kx)+\bg_n(\kx,\bo),
\quad n=0,1,\dots,\, \kx\in\z2.
\endd
\bigskip
Now we are ready to prove our
\begin{lemma} For any symmetric, recurrent aperiodic walk in
two dimension we have

\begg \pe(\bt_a\geq k) \leq 2\frac{g(a)}{g(k-1)},
\quad a=0,1,\dots,\, k=1,2,\dots \endd
\end{lemma}

\bigskip
\nod
{\bf Proof:} We start with the following simple observation: for $k\geq
1$, $a\geq 0$

\begg \pe(\bt_a\geq k)=\sum_{\kx \in \z2} \pe(\bt_a\geq
k|\bs_a=\kx) \pe(\bs_a=\kx)=\sum_{\kx\in
\z2}P_a(\bo,\kx)\pe_{\kx}(\bt\geq k). \label{lem20}
\endd

Now (\ref{sp11})  and the symmetry of the walk implies
that
\begg a(\kx+\ky)\leq a(\kx)+a(\ky). \label{lem22}
\endd

Moreover, denoting $\kx+\ky=\kz$ we easily get from (\ref{lem22})
that
\begg
 a(\kz)\leq a(\kx)+a(\kz-\kx)=a(\kx)+a(\kx-\kz),
\endd
thus
\begg
a(\kz)-a(\kx-\kz)\leq a(\kx),
\endd
which in turn, combined with (\ref{sp11}) implies that \begg
g(\kx,\ky )\leq 2\,a(\kx). \label{lem23}
\endd
Combining (\ref{sp6}), (\ref{sp7}) and (\ref{lem23}) we conclude
that for $\kx \neq\bo$ \begg \pe_{\kx}(\bt>n)\leq 2a(\kx)
\pe(\bt>n+1),\quad n=0,1,\dots \label{lem24}
\endd
Now using (\ref{lem20}),  (\ref{lem24}), Fact 4 and  the simple
observation that $a(\bo)=0$, we get that
\begin{eqnarray}
\pe(\bt_a\geq k)&\leq& 2 \pe( \bt\geq k+1)\sum_{\kx\in
\z2-\{\bo\}}P_{a}(\bo,\kx)a(\kx)+P_{a}(\bo,\bo)\pe(\bt\geq
k)\nonumber\\
&\leq& 2\pe(\bt \geq k)
\le\sum_{\kx\in\z2}P_{a}(\bo,\kx)a(\kx)+P_{a}(\bo,\bo)\ri\nonumber\\
&\leq&
2\pe(\bt \geq k) (\bg_{a-1}(\bo,\bo) +P_{a}(\bo,\bo))\leq
2\pe(\bt\geq k)g(a).\label{lem25}
\end{eqnarray}

\nod
To estimate $\pe(\bt> k)$ we use an argument essentially from
Erd\H os and Taylor \cite{ET60}. Partitioning according to the last
return to zero we get
\begg
\sum_{j=0}^k \pe(\bt>k-j)P_j(\bo,\bo)=1
\endd
implying that
 \begg \pe (\bt>k)\sum_{j=0}^k P_j(\bo,\bo)\leq 1,
\endd
hence \begg \pe (\bt>k)\leq 1/g(k). \label{gn1}
\endd
Now (\ref{lem25}) and (\ref{gn1}) imply our lemma. $\Box$

\medskip
\noindent {\bf Fact 5.} { \it If the walk is recurrent and
aperiodic with finite second moment, then we  have as}
$n\to\infty$ \begg g(n)\sim c_2 \log n \qquad {\rm for\,\,\,}
d=2,\endd
 \begg g(n)\sim c_1 \sqrt{n}
\qquad {\rm for\,\,\,} d=1.\endd

\bigskip
 The case $d=2$ is well-known for strongly
aperiodic walk (see \cite{S76}, p. 75) and \cite{DPRZ01} how to
weaken this condition for the aperiodic case). The case  $d=1$  is
well-known (see e.g. \cite{S76}, p. 381).

 Let $\{\bs_n\}$ be a symmetric recurrent  aperiodic
walk in two dimension. Consider the following problem. Let
$\bt_a$ be defined by (\ref{sp1}). At $a+\bt_a$ the walk is at $\bo$.
Now after another $a$ steps  we wait again until the walk
returns to $\bo.$ Keep repeating this procedure, we would like to estimate the
number of such returns within $n$ steps.  This problem was
considered in \cite{CSF83} in one dimension. Here we repeat
essentially the same argument  and spell out it in the  $d=2$ case
with the appropriate modifications, using Lemma 2.2. To formulate
this problem precisely, let
 \begg
\zeta_1(a):=a, \quad\alpha_k(a):=\bt_{\zeta_k(a)},
\quad\zeta_{k+1}(a):=\sum_{i=1}^{k}\alpha_i(a)+(k+1)a,\quad
k=1,2,\dots
\endd
Then $\alpha_k(a),\, k=1,2,\dots$ are i.i.d. random variables having the
same distribution as $\bt_a$.

\begin{lemma} For a symmetric recurrent aperiodic random walk in $\z2$
we have for $a>1,$ $u>1$, $k\geq 1$
\begg \pe\left(\sum_{i=1}^k \alpha_i(a)\geq u\right)\leq
C\,k \frac{\log a}{\log u}.
\endd
\end{lemma}
{\bf Proof:} By Lemma 2.2  and Fact 5

\begg \pe(\bt_a\geq k) \leq 2\frac{g(a)}{g(k-1)} \leq C
\frac{\log a}{\log(k-1)}\leq C \frac{\log a}{\log k}.
\endd

Let
\begg
^u\alpha_k(a)=\left\{\begin{array}{ll}\alpha_k(a)
\quad &{\rm if \quad}\alpha_k(a)\leq u \\0
\quad & {\rm if\quad}\alpha_k(a)> u \\
\end{array}\right.
\label{main3}
\endd

First observe, that
 \begin{eqnarray} \be(^u\alpha_i(a))\leq
  \sum_{j=0}^{u}
\pe(\bt_a\geq j) \leq  2+C \sum_{j=2}^u\frac{\log a}{\log
j}\leq
 C u\frac{\log a}{\log u}.
\end{eqnarray}

 Define the event
 \begg A=\bigcap_{j=1}^k\{\alpha_j(a)<u\}.
\endd

\begin{eqnarray}
\pe\left(\sum_{i=1}^k \alpha_i(a)\geq u\right)&\leq&
\pe\left(\sum_{i=1}^k \alpha_i(a)\geq u, A\right)+\pe(\bar
{A}) \leq \pe\left(\sum_{i=1}^k\, ^u\alpha_i(a)\geq
u\right)\nonumber\\
&+& k \pe(\alpha_i(a)\geq u)\leq
\frac{k\be(^u\alpha_i(a))}{u}+k C\frac{\log a}{\log u}\leq C
k \frac{\log a}{\log u}.
\end{eqnarray}
$\Box$

Let $a_t>1$ be an integer valued function of $t$ and let $\nu_t$ be the
largest integer $N$ for which
\begg
\sum_{i=1}^N \alpha_i(a_t)+(N+1)
a_t\leq t.
\endd
\begin{lemma} Under the conditions of Lemma 2.3 for $f(t)>0,$
$(f(t)+2)a_t<t-1$ we have

 \begg \pe (\nu_t\leq f(t))\leq
 C f(t)\, \frac {\log a_t}{\log(t-(f(t)+2)a_t)}.
\endd
\end{lemma}
{\bf Proof:} By Lemma 2.3
\begin{eqnarray}
\pe (\nu_t\leq f(t))&\leq& \pe (\nu_t\leq [f(t)]+1)\leq
\pe\left(\sum_{i=1}^{[f(t)]+1}\alpha_i(a_t)>t-([f(t)]+2)a_t)\right)
\nonumber\\
&\leq&  C(f(t)+1)\frac {\log
a_t}{\log(t-(f(t)+2)a_t)},
\label{lem241}
\end{eqnarray}
implying our statement. $\Box$

As mentioned above, the quantities defined in this section have analogues
in one dimension. Also, there are corresponding one-dimensional analogues
of results quoted for two-dimensional case. We recall the
one-dimensional versions of Lemmas 2.3 and 2.4.

\bigskip

\noindent {\bf Fact 6.} (\cite{CSF83}) {\it Consider a symmetric
aperiodic random walk $\{S_n\}$ on ${\cal Z}_1$ with finite
variance and define $\al_i(a)$ and $\nu_t$ exactly as before for
$\{S_n\}$. Then for} $a>1$, $u>1$, $k\geq 1$

\begg \pe\left(\sum_{i=1}^k \alpha_i(a)\geq u\right)\leq
C\,k \sqrt{\frac{a}{u}}.
\endd
{ \it Furthermore, if $f(t)>0,$ $a_t>0,$ $(f(t)+2)a_t< t$ we have}

 \begg \pe (\nu_t\leq f(t))\leq
C f(t)  \frac {\sqrt{ a_t}}{\sqrt{t-(f(t)+2)a_t}}.
\endd
We will need the following upper tail estimates  essentially from
Erd\H{o}s and Taylor \cite{ET60}.

\bigskip
\noindent {\bf Fact 7.} (\cite{ET60} or \cite{DPRZ01}) {\it For the simple
symmetric random walk on the plane we have for any  $\alpha>0$ and}
$0<\delta<1$, $n\geq 1$

\begg \pe(\xi^{(2)}(\bo,n)\geq \al(\log n)^2)< n^{-(1-\delta)\pi
\al}.
\endd

 \bigskip
 \nod
{\bf Fact 8.} (\cite{ET60}, (3.6))
{\it For the simple symmetric random walk on the plane we have for any}
$\alpha>0$, $\delta>0$, $n\geq 1$

 \begg
\pe(\xi^{(2)}(\bo,n)\geq \al(\log n)^2)> n^{-(1+\delta)\pi \al}.
\endd
\bigskip
 \nod
{\bf Fact 9.} (\cite{CSF83})
{\it Let} $\{S_n\}$ {\it be a one-dimensional symmetric aperiodic random
walk on} ${\cal Z}_1$ {\it with}
  $\sigma^2=E(X_1^2)<\oo$ {\it and suppose that} $x_n\to
  \oo,\,\,x_n/n^{1/2}\to 0$ {\it as} $\noo.$ {\it Let} $\xi(0,n)$ {\it be
its local time at zero. Then for any} $\ep>0$ {\it and large enough} $n$
\begg
C_1\exp\le-(1+\ep)\frac{x_n^2\sigma^2}{2}\ri\leq \pe(\xi(0,n)\geq
x_n\, n^{1/2})\leq C_2 \exp\le-(1-\ep)\frac{x_n^2\sigma^2}{2}\ri.
\endd
{\bf Fact 10.} (\cite{CSRR98}, (2.1) and \cite{MR94}, Lemma 2.5.)
{\it For a symmetric aperiodic random walk with finite variance in
$\z2$ we have for any $x>0$}

\begg \pe\le \xi^{(2)}(\bo, n)\geq x\log n \ri\leq
\exp(-c x)
\endd
with some constant $c>0$.

\renewcommand{\thesection}{\arabic{section}.}
\section{Proofs of the two dimensional results.}
\renewcommand{\thesection}{\arabic{section}}
\setcounter{equation}{0} \setcounter{theorem}{0}
\setcounter{lemma}{0}

{\bf Proof of Theorem 1.1.} Consider the line  $L=L(a_1,a_2)$ with
$a_1x_1+a_2x_2=0$, where $a_1$ and $a_2$ are integers, not both of them
zero. Without loss of generality we may assume that $a_1$ and $a_2$ are
relatively prime.

For the two-dimensional random walk define a one-dimensional walk with
steps

$$ Y_i=a_j\quad {\rm if \quad} \bx_i=\bee_j,\qquad i=1,2,...\quad
j=1,2.
$$

$$ Y_i=-a_j\quad {\rm if \quad} \bx_i=-\bee_j,\qquad i=1,2,...\quad
j=1,2.
$$
Then $Z_n=\sum_{i=1}^n Y_i$ is an aperiodic one-dimensional
symmetric random walk with $Z_n=0$ if and only if $\bs_n\in L(a_1,a_2).$ Thus
denoting by $V^L(n)=\#\{i: 1\leq i\leq n, \,\, \bs_i\in L(a_1,a_2)\}, $
the number of visits of $\bs_n$ on $L(a_1,a_2),$ we have

\begg V^L(n)=\xi^{Z}(0,n), \label{111} \endd
where $\xi^{Z}(\cdot,n)$ is the  local time of $\{Z_n\}.$

To get the upper bound in our theorem, denote by $D(n)$ the
set of lattice points on $L$ which are visited up to $n$  by
$\{\bs_i\}$, and denote by $|D(n)|$ the number of points in $D(n).$
Select a subsequence $n_j=[e^j]$, $j=1,2,...$ and let for
any $ \kx \in \z2$
$$C_{\kx}^j=\{ \xi^{(2)}(\kx,n_{j+1})>\lambda(\log n_j)^2\}.$$

\nod Then using  (\ref{111}), Facts 7 and 9 we conclude that for
any $\lambda>0$, $\delta>0$, $\varepsilon>0$ and $j$ large enough

\begin{eqnarray}
&&\pe(\xi^{(2)}_L(n_{j+1})>\lambda(\log n_j )^2)\nonumber\\
&\leq& \pe\le \bigcup_{\kx
\in D(n_{j+1})}C_{\kx}^j,\,\,\, |D(n_{j+1})|\leq \sqrt{ n_{j+1}\log
n_{j+1} }\ri+  \pe \le |D(n_{j+1})|> \sqrt{n_{j+1}\log
n_{j+1}}\ri \nonumber\\
&\leq&
  \pe\le \bigcup_{\kx \in D(n_{j+1})}C_{\kx}^j,\,\,\,
|D(n_{j+1})|\leq \sqrt{n_{j+1}\log n_{j+1}}\ri +  \pe\le
\xi^{Z}(0,n_{j+1})> \sqrt{n_{j+1}\log n_{j+1}}\ri\nonumber\\
 &\leq& \sqrt{n_{j+1}\log n_{j+1}}\,\,\,\pe(\xi^{(2)}(\bo,n_{j+1})>
\lambda(\log n_j)^2) +C_2 \exp\le-(j+1) \frac{(1-\ep) \sigma^2}{2}\ri
\nonumber\\
&\leq&\sqrt{n_{j+1}\log
n_{j+1}}\,\exp\le-\lambda(1-\delta)\pi\frac{j^2}{j+1}\ri
+ C_2 \exp\le-j \frac{(1-\ep) \sigma^2}{2}\ri\nonumber\\
&&\leq\exp \le-(1-\delta)^2\pi \,j\,\lambda+j(1/2+\ep)\ri +
C_2 \exp\le-j \frac{(1-\ep) \sigma^2}{2}\ri,
\label{112}
\end{eqnarray}
where we used  in the second inequality above
that $|D(n)|\leq V^L(n) .$ $\sigma^2$ is the variance of $Y_i$ and hence
depends only on $a_1$ and $a_2.$ It is easy to see that we can choose
$\ep>0$, $\delta>0$ for which the last line of (\ref{112}) is
summable in $j$ whenever $\lambda>1/ (2\pi)$, which in turn, using
Borel-Cantelli lemma and the usual monotonicity argument implies
\begg
\limsup_{\noo}\frac{\xi^{(2)}_L(n)}{(\log n)^2}\leq \frac{1}{2\pi},
\endd
so we have the upper half of the theorem.

\bigskip
To get the lower bound, we will essentially follow Erd\H{o}s and
Taylor's argument  with some modification. We consider the walk
$\{\bs_i\}$,  wait $[n^{\al}]$ steps and  observe $\xi^{(2)}(\bo,[n^{\al}]).$
The number $0<\al<1$ will be choosen later. At time $[n^{\al}]$ the walk
is somewhere on the plane, and we wait until its first return to
$L=L(a_1,a_2).$ When it returns to $L$, we observe the local time of the
hitting point of $L$ for a time interval $[n^{\al}],$ and then wait again
for the walk to return to $L. $ We keep repeating this procedure for a
total time of $n$ steps. This construction ensures that the local times of
these hitting points over a time interval $[n^\al]$ are  i.i.d. random
variables having the same distribution as $\xi^{(2)}(\bo,[n^{\al}])$.
Combining this observation with our Fact 6, would produce  our lower bound
in the theorem. In what follows we work out the above outlined ideas with
the added complexity of working with subsequences as before. Let
$n_j=j^{\beta},\quad j=1,2,\dots$ with integer $\beta$ to be choosen
later. Now fix $j\geq 1$ and define for $k=1,2,\dots$

\begg
\eta_0^j:=0,\qquad \eta_k^j:=\min\{i>\eta_{k-1}^j+[n^\al_j]:\,
\bs_i\in L\},\qquad \by_k^j:=\bs_{\eta_k^j}.
\endd

Furthermore, let $\nu_j$ be the largest integer $N$ for which
\begg
\eta_N^j+[n_j^\al] \leq n_j.
\endd

Put $f(n_j)=[n_j^{(1-\al)(1-\eta)/2}]$ for any $0<\eta<1$. Since $\bs_i\in
L$ if and only if $Z_i=0$, we can apply Fact 6 to get

\begg
\pe \le\nu_j<f(n_j)  \ri\leq C n_j^{(\al-1)\eta/2}.
\label{arr1}
\endd

Introduce further the events

\begg A_k^j=\{
\xi^{(2)}(\by_{k}^j,\eta_k^j+[n_j^\al])-\xi^{(2)}(\by_{k}^j,\eta_{k}^j)<K
(\log n_{j+1})^2\}.
\endd
For fixed $j\,$ the events $A_k^j$ are independent in $k$ and
having the same probability as $A_1^j.$ Using the above notations and
(\ref{arr1})

\begin{eqnarray}
\pe \le \xi_L^{(2)}(n_j)<K (\log n_{j+1})^2 \ri\leq \pe\le
\bigcap_{k=1}^{\nu_j} A_k \ri \leq C
n_j^{(\al-1)\eta/2}+\le\pe(A_1^j)\ri^{f(n_j)}.
\end{eqnarray}
 Now by Fact 8 we have
 \begin{eqnarray}
\pe(A_1^j)&=&\pe(\xi^{(2)}(\bo,[n_j^\al])<K(\log n_{j+1})^2)
\nonumber\\&\leq&1-\exp\left\{-\frac{K\pi}{\al}(1+\delta)(\beta\log j)
\le\frac{\log(j+1)}{\log j}\ri^2\right\}.
 \end{eqnarray}
Consequently, we have for $j$ big enough
\begin{eqnarray}
&&\pe \le \xi_L^{(2)}(n_j)<K (\log n_{j+1})^2 \ri\nonumber\\
&\leq&Cn_j^{(\al-1)\eta/2}+
\le1-\exp\left\{-\frac{K\pi}{\al}(1+\delta)^2(\beta\log
j) \right\}\ri^{(n_j)^{\frac{(1-\al)(1-\eta)}{2}}}\nonumber\\
&\leq&Cj^{\beta(\al-1)\eta/2}+\le
1-\frac{1}{j^{\frac{K\pi}{\al}\beta(1+\delta)^2}} \ri
^{j^{\frac{\beta(1-\al)(1-\eta)}{2}}}\nonumber\\
&\leq& Cj^{\beta(\al-1)\eta/2}+
C\exp\{-j^{\frac{\beta(1-\al)(1-\eta)}{2}-\frac{K\pi}{\al}\beta(1+\delta)^
2}\}.
\label{vegre}
\end{eqnarray}

\nod
For given $0<\al<1$ and $\eta>0$ select an integer $\beta$ such that
$$
\beta(\al -1)\eta/2<-1
$$
so that the first term in (\ref{vegre}) is summable in $j$.
The second term will be summable in $j$ whenever
\begg
\frac{(1-\al)(1-\eta)}{2}>\frac {K\pi}{\al}(1+\delta)^2.\endd
On choosing $\alpha=1/2$, $\delta>0$ and $\eta>0$  small, we conclude that
(\ref{vegre}) is summable in $j$ if
\begg K <
\frac{1}{8\pi}.
\endd
Using again Borel-Cantelli lemma and the usual monotonicity argument,
we get

\begg \limsup_{\noo}\frac{\xi^{(2)}_L(n)}{(\log n)^2}\geq \frac{1}{8\pi},
\endd
proving our theorem. $\Box$

\bigskip
\nod  {\bf Proof of Theorem 1.2.}   First observe that the
condition $\al \leq 1/2$ is not a real restriction as by the LIL
(Fact 1) with probability 1 the walk in the time interval $[0,n]$ remains
in $B(r_n)$ for $r_n=n^{\al},\quad \al> 1/2$, $n$ large, thus
$\xi^{(2)}_{B(r_n)}(n)=\xi^{(2)}(n)$ eventually with probability 1.

In case ${\al \leq 1/2}$, using the LIL again, during the first
$[r_n^{2-\ep}]$ steps the walk remains in $B(r_n)$ with probability 1 for
any $\ep >0$ and large $n$. Hence by Theorem C for any $\ep>0$ and large
enough $n$,
\begin{eqnarray}
\xi^{(2)}_{B(r_n)}(n)&\geq& \xi^{(2)}(r_n^{2-\ep})\geq
\frac{1}{\pi} \le \log r_n^{2-\ep} \ri^2(1-\ep) \nonumber
\\ &=&\frac{1}{\pi} \le \log
n^{\al(2-\ep)}\ri^2(1-\ep)=\frac{\al^2(2-\ep)^2}{\pi} (\log
n)^2(1-\ep) \has
\end{eqnarray}

\nod Let $\ep \to 0$  to get

\begg  \liminf_{\noo} \frac{\xi^{(2)}_{B(r_n)}(n)}{(\log n)^2}
\geq\frac{4 \al^2}{\pi} \has
\endd

The upper bound is proved in Dembo {\it et al.} \cite{DPRZ01}.

Thus we have our theorem. $\Box$

\nod
{\bf Remark.} Observe that in case $\al=1/2$  the $\liminf$ and
$\limsup$ coincide.

\bigskip

\nod  {\bf Proof of Theorem 1.3.}   Let $r_n= \exp\{(\log n
)^\beta\}.$  To get the lower half of the theorem, by  the same
argument as in  Theorem 1.2,  in the first $[r_n^{2-\ep}]$  steps the
walk remains in $B(r_n)$ with probability 1, thus by Theorem C
for large $n$ we have

\begin{eqnarray}
\xi^{(2)}_{B(r_n)}(n)&\geq& \xi^{(2)}([r_n^{2-\ep}])\geq
\frac{1}{\pi} \le \log r_n^{2-\ep} \ri^2(1-\ep)\nonumber\\
&=&
 \frac{1}{\pi}  (\log n)^{2 \beta}(2-\ep)^2(1-\ep) \has
\end{eqnarray}

\nod Sending now $\ep \to 0$ we get

\begg  \liminf_{\noo} \frac{\xi^{(2)}_{B(r_n)}(n)}{(\log n)^{2\beta}}
\geq\frac{4 }{\pi} \has
\endd

To get the other half of the theorem, consider first a disc
$A_n:=B(\exp\{ n^{\gamma}\})$ where $1/2\leq \gamma < 1$. Recall the
definition of $\rho_n$ in (\ref{rho}). First we show that for arbitrary
small $\delta>0$

\begg
\max_{\kx \in A_n} \xi^{(2)}(\kx, \rho_n) \leq n^{2\gamma+\delta} \has
\label{130}
\endd

\nod
Using the definition of $Y_n(\kx)$ in Fact 3 and  Lemma 2.1  we have

\begin{eqnarray}
\pe(\max_{\kx \in A_n}\, \xi^{(2)}(\kx, \rho_n)> u \, n) &\leq&
\pe(\max_{\kx \in A_n} \, \sum_{k=1}^n Y_k(\kx)>u\,n )\nonumber\\
&\leq&  C_1 \exp\{2 n^{\gamma} \}\,\max_{\kx \in A_n}\, \pe(\sum_{k=1}^n
Y_k(\kx)>u\,n )\nonumber \\
&\leq&  C_1 \exp\{2 n^{\gamma}
\}\,\max_{\kx \in A_n}\, \exp\{ n p(\kx)(1-u/2) \} \nonumber\\
&\leq&
C_1 \exp\{2 n^{\gamma} \} \exp\{C_2 n^{1-\gamma}(1-u/2) \},
\label{131}
\end{eqnarray}
where in the last inequality we used that for all  $\kx \in
A_n\,\,$ $p(\kx)\geq C_2 n^{-\gamma}  $ by Fact 2. Selecting now
$u=n^{\theta}$ we get from (\ref{131}) that for $n$ big enough
\begg \pe(\max_{\kx \in A_n} \xi^{(2)}(\kx, \rho_n)>
n^{1+\theta}) \leq C_1 \exp\{2 n^{\gamma} \} \exp\{-C_3
n^{1+\theta-\gamma}\}. \label{132}
\endd
\nod The probabilities in (\ref{132}) can be summed up for $n$ if
$1+\theta-\gamma> \gamma.$ Thus for any  $1/2\leq\gamma<1$ and for an
arbitrary small $\delta>0$ we might select $\theta >0$  such that

\begg
1+\theta= 2\gamma+\delta
\endd
to get by Borel-Cantelli lemma that (\ref{130}) holds.

Applying  now (\ref{130}) with  $n=\xi^{(2)}(\bo, k) $ we get that

\begin{eqnarray}
\max_{\kx \in A_{\xi^{(2)}(\bo, k)}} \xi^{(2)}(\kx, k)&\leq& \max_{\kx \in
A_{\xi^{(2)}(\bo, k)+1}}
\xi^{(2)}(\kx, \rho_{\xi^{(2)}(\bo, k)+1})\nonumber \\
&\leq &  \le \xi^{(2)}(\bo, k)+1 \ri ^{2\gamma+\delta}
\leq (\log k)^{2\gamma+\delta'} \has
\label{133}
\end{eqnarray}
for any $\delta'>\delta$ by Theorem A. On the
other hand, (\ref{133}) and Theorem A imply that
for any $0<\gamma'<\gamma$ and $k$ large enough

\begin{eqnarray}
\max_{\kx \in B(\exp\{ (\log k)^{\gamma'}\})} \xi^{(2)}(\kx, k)\leq
\max_{\kx \in A_{\xi^{(2)}(\bo, k)}} \xi^{(2)}(\kx, k)
\leq(\log k)^{2\gamma+\delta'} \has
\end{eqnarray}
which is equivalent to our statement.
$\Box$

\bigskip
\nod {\bf Proof of Corollary 1.1.} Clearly the  upper bound for
$\xi^{(2)}_{B(r_n)}(n)$ or $\xi^{(2)}_L(n)$ holds for
$\xi^{(2)}_{B(r_n)}\cap L$ as well. On the other hand, to get the
lower bounds, observe that with probability one, in the first
$[r_n^{2-\ep}]$ steps the walk remains in $B(r_n)$ for large $n$.
Consequently \begg \xi^{(2)}_L(r_n^{2-\ep})=\xi^{(2)}_{L \cap
B(r_n)}(r_n^{2-\ep})\leq \xi^{(2)}_{L\cap B(r_n)}(n). \label{cor}
\endd
But (\ref{cor}) and Theorem 1.1 imply our statements. $\Box$

\renewcommand{\thesection}{\arabic{section}.}
\section{Proofs of the higher dimensional results.}
\renewcommand{\thesection}{\arabic{section}}
\setcounter{equation}{0} \setcounter{theorem}{0}
\setcounter{lemma}{0}

{\bf Proof of Theorem 1.4}. The proof of this theorem
is very similar to the proof  of Theorem 1.1. Consider the subspace

$$ S_{d-1}=\{\kx\in\zd:\, a_1x_1+ a_2x_2+...+a_dx_d  =0   \}
$$
with integer coefficients  $a_1,\,a_2,\,...a_d,\,$ not all of them zero.
Without loss of generality we may assume that the largest common divisor
of $(a_1,\dots,a_d)$ is equal to 1.

Define a one-dimensional random walk with the following steps:
 $$
Y_i=a_\ell \quad {\rm if}\quad \bx_i= \bee_{\ell},\qquad \ell=1,2,
...d, \qquad i=1,2,...
$$
 \begg
Y_i=-a_\ell \quad {\rm if}\quad \bx_i= -\bee_{\ell},\qquad
\ell=1,2, ...d, \qquad i=1,2,...
\endd
Note that these values are not necessarily distinct. In that case
we sum up the probabilities $1/(2d)$ according to their
multiplicity.

$Z_n=\sum_{i=1}^n Y_i,\, n=1,2,\dots$ is an aperiodic one-dimensional
symmetric random walk with $Z_n=0$ if and only if $\bs_n\in S_{d-1}.$
 Thus denoting by $V^{S_{d-1}}(n)=\#\{i: 1\leq i\leq n, \,\, \bs_i
\in S_{d-1} \}, $ the number of visits of $\{\bs_i\}$ up to time $n$ in
$S_{d-1},$ we have

\begg V^{S_{d-1}}(n)=\xi^{Z}(0,n), \label{411} \endd
where
$\xi^{Z}(0,n)$ is the local time at zero up to time $n$ of the
random walk $\{Z_i\}$.

To get the upper bound in our theorem, denote by $D(n)$ the  set
of lattice points on $S_{d-1}$ which are visited by $\{\bs_i\}$ up to
time $n$ and denote $|D(n)|$ the number of points in $D(n)$.
Select a subsequence $n_j=[e^j], \quad j=1,2,\dots$ and define the
events for any $ \kx\in\zd $
$$C_{\kx}^j=\{ \xi^{(d)}(\kx,n_{j+1})>\lambda (\log n_j)\}.$$

\nod Then using the exact distribution (\ref{geo}), and
 Fact 9  we conclude, that

\begin{eqnarray}
&&\pe(\xi^{(d)}_{S_{d-1}}(n_{j+1})>\lambda\log n_j )
\nonumber\\
&\leq& \pe\le
\bigcup_{\kx \in D(n_{j+1})}C_{\kx}^j,\,\,\, |D(n_{j+1})|\leq
\sqrt{n_{j+1}\log n_{j+1} }\ri + \pe \le |D(n_{j+1})|>
\sqrt{n_{j+1}\log n_{j+1}}\ri\nonumber\\
 &\leq&
  \pe\le \bigcup_{\kx \in D(n_{j+1})}C_{\kx}^j,\,\,\,
|D(n_{j+1})|\leq \sqrt{n_{j+1}\log n_{j+1}}\ri +  \pe\le
\xi^{Z}(0,n_{j+1})> \sqrt{n_{j+1}\log n_{j+1}}\ri \nonumber\\
&\leq& \sqrt{n_{j+1}\log n_{j+1}}\,\,\,\pe(\xi^{(d)}(\bo,n_{j+1})>\la \log
n_j) + C_2 \exp\le-(j+1) \frac{(1-\ep) \sigma^2}{2}\ri \nonumber\\
&\leq& \sqrt{n_{j+1}\log n_{j+1}}\,\pe(\xi^{(d)}(\bo,\infty)>\la\log n_j)
+C_2 \exp\le-(j+1) \frac{(1-\ep) \sigma^2}{2}\ri \nonumber\\
&\leq& \sqrt{n_{j+1}\log n_{j+1}}\exp \le -\frac{\la}{\ld}j\ri +  C_2\exp
\le-(j+1) \frac{(1-\ep) \sigma^2}{2}\ri \nonumber\\
&&\quad\leq
C\sqrt{j}\exp\le j\le\frac{1}{2}-\frac{\la}{\ld}\ri\ri +C_2\exp
\le-(j+1) \frac{(1-\ep) \sigma^2}{2}\ri
 \label{412}
\end{eqnarray}
for  $j$ big enough, where we used  in the second inequality above
that $|D(n)|\leq V^{S_{d-1}}(n).$ The $\sigma$ above depends only
on $a_1, a_2,...a_d.$
 It is easy to see that the  last line of (\ref{412}) is
summable in $j$ whenever $\la>\ld/2$, which in turn, using
Borel-Cantelli lemma and the usual monotonicity argument implies
\begg
\limsup_{\noo}\frac{\xi^{(2)}_{S_{d-1}}(n)}{\log n}\leq
\frac{\ld}{2} \has, \label{sup}
\endd
 so we have the upper half of the theorem.

To get the lower bound, we will again follow Erd\H{o}s and
Taylor's argument  with the appropriate modification. We consider the
walk $\{\bs_n\}$,  wait $[(\log n)^2]$ steps. After that many steps we
wait
until the walk arrives back to $S_{d-1}$. Then wait again $[(\log n)^2]$
steps, and repeat this procedure over and over again. The probability that
an arrival point to $S_{d-1}$ will be visited by the walk again in the
next $[\log n]$ steps is by (\ref{gam2})
$$
1-\gamma_d+O\left(\frac{1}{(\log n)^{1/2}}\right)
$$
and the probability that it will be visited  at least $[\la \log
n]$ times (with some  $\lambda<1$) within $[(\log n)^2]$ steps is
greater than \begg \le 1-\gamma_d+O\left(\frac{1}{(\log
n)^{1/2}}\right)\ri^{[\la \log n]}. \label{gam100}
\endd
These ideas will be combined with our Fact 6, and applied for subsequences,
just like in the proof of Theorem 1.1.

Let  $n_j=j^{\beta},\quad j=1,2,\dots$ with some integer $\beta$. Define
\begg
\eta_0^j:=0,\qquad \eta_k^j:=\min\{i>\eta_{k-1}^j+[(\log n_j)^2]:\,
\bs_i\in S_{d-1}\},\qquad \by_k^j:=\bs_{\eta_k^j}.
\endd

Furthermore, let $\nu_j$ be the largest integer $N$ for which
\begg
\eta_N^j+ [(\log n_j)^2] \leq n_j.
\endd

Selecting for any   $0<\eta<1, \,\,$
$f(n_j)=\left[\le \frac{n_j}{( \log n_j)^{2}}\ri ^{(1-\eta)/2}\right]$, it
is easy to see that from Fact 6 we get

\begg \pe \le\nu_{j}<f(n_j)  \ri\leq C \le \frac{( \log n_j)^{2}}{n_j}\ri
^{\eta/2}.
\label{arr11}
\endd
\nod
Introduce further the events

\begg A_k^j=\{
\xi^{(d)}(\by_{k}^j,\eta_k^j+[(\log
n_j)^2])-\xi^{(d)}(\by_{k}^j,\eta_{k}^j)<
\lambda \log n_{j+1}\}.
\endd

\nod
For fixed $j\,$ the events $A_k^j$ are independent in $k$ and
 having the same probability as $A_1^j.$   Using (\ref{arr11}) we get that

\begin{eqnarray}
&&\pe \le \xi^{(d)}_{S_{d-1}}(n_j)<\lambda \log n_{j+1} \ri\leq \pe\le
\bigcap_{k=1}^{\nu_{j}} A_k^j \ri\nonumber\\
 &\leq& C \le \frac{( \log
n_j)^{2}}{n_j}\ri
^{\eta/2}
+\le\pe(A_1^j)\ri^{f(n_j)}.
\end{eqnarray}
\nod
Now by (\ref{gam100}) we have
\begin{eqnarray}
&&\pe(A_1^j)=\pe(\xi^{(d)}(\bo,[(\log n_j)^2])<\lambda(\log n_{j+1}))
\nonumber\\
&<&1-\le 1-\gamma_d+O\left(\frac{1}{(\log n_j)^{1/2}}\right)\ri^{\la \log
n_{j+1}} =1-n_{j+1}^{-\frac{\lambda}{\ld}\left(1+O\left(\frac{1}{(\log
n_j)^{1/2}}\right)\right)}.
\end{eqnarray}
Consequently,

 \begin{eqnarray}
&&\pe \le \xi^{(d)}(n_j)<\lambda (\log n_{j+1}) \ri\nonumber\\
&\leq&C \le\frac{\beta^2(\log j)^2}{j^{\beta}}\ri^{\eta/2}
 +\le1-n_{j+1}^{-\frac{\lambda}{\ld}(1+O(\frac{1}{(\log
n_j)^{1/2}}))}
 \ri^{\left[\le \frac{n_j}{( \log n_j)^{2}}\ri ^{\frac{1-\eta}{2}}\right
 ]}\nonumber\\
&\leq&
 j^{-\frac{\beta \eta}{3}}+ \le 1-(j+1)^{-
\beta\frac{\lambda}{\ld}(1+O(\frac{1}{(\log
j)^{1/2}}))}\ri^{j^{\frac{\beta(1-2\eta)}{2}}}\nonumber\\
&\leq&
j^{-\frac{\beta \eta}{3}}+\exp\le - \frac{1}{2}j^
{\beta\le\frac{1}{2}-\eta- \frac{\lambda}{\ld} +O (\frac{1}{(\log
j)^{1/2}})\ri} \ri \label{finally}
\label{4upper}
\end{eqnarray}
for $0<\eta<1/2$ and $j$ big enough.
 It is easy to see that in (\ref{finally}) selecting  $\eta>0$ as small as
necessary
 the second term is summable in $j$ if $\la<\ld/2.$ On the other hand, one
can select a $\beta$ big
  enough (depending on $\eta$) that the first term is summable in $j.$
Borel-Cantelli lemma  and
  the usual monotonicity argument results in
\begg \liminf_{\noo}\frac{\xi^{(2)}_{S_{d-1}}(n)}{\log n}\geq
\frac{\ld}{2} \has
\endd
which combined with (\ref{sup}) implies our theorem.
$\Box$

\bigskip\noindent
{\bf Proof of Theorem 1.5}. Recall that
$$
S_{d-2}=\{\kx\in\zd:\, a_1x_1+ a_2x_2+...+a_dx_d  =0, \quad b_1x_1+
b_2x_2+...+b_dx_d=0\}
$$
with integers $a_i, b_j$. Without loss of generality, we may assume that
the largest common divisor of both $(a_1,\dots,a_d)$ and $(b_1,\dots,b_d)$
is equal to 1 and $(a_1,\dots,a_d)\neq(b_1,\dots,b_d)$.

Define a two-dimensional random walk associated
with our $d$-dimensional random walk with steps
$$
\widetilde Y_i=(a_r,b_r)\quad \textrm{if}\quad \bx_i=\bee_r,\quad
r=1,2,\dots, d,
$$
$$
\widetilde Y_i=-(a_r,b_r)\quad \textrm{if}\quad \bx_i=-\bee_r,\quad
r=1,2,\dots,d.
$$
Note that, as in the previous theorem, these values are not
necessarily distinct. In that case we sum up the probabilities
$1/(2d)$ according to their multiplicity.

Then
$$
\widetilde Z_n=\widetilde Y_1+\dots+\widetilde Y_n,\quad
n=1,2,\dots
$$
is a two-dimensional symmetric recurrent (possibly periodic) random walk
in $\z2$ with finite variance with the property $\bs_n\in
S_{d-2}$ if and only if $\widetilde Z_n=0$. Let
$$
V^{S_{d-2}}(n)=\#\{i:\, 1\leq i\leq n,\, \bs_i\in S_{d-2}\}
$$
and

$$
\xi^{\widetilde Z}(\bo, n)=\#\{i:\, 1\leq i\leq n,\, \widetilde Z_i=\bo\},
$$
the local time of $\{\widetilde Z\}$ at $\bo$ up to time $n$.

Then we have $V^{S_{d-2}}(n)=\xi^{\widetilde Z}(\bo, n)$.

Now we define an aperiodic random walk with the same property. If
$\{\widetilde Z_n\}$ is aperiodic, then let $\widehat Z_n=
\widetilde Z_n,\, n=0,1,\dots$ In the case when $\{\widetilde Z_n\}$
happens to be periodic, then following Spitzer \cite{S76}, pp. 65-66
construct an aperiodic walk in $\z2$ as follows. The steps $(a_r, b_r),\,
r=1,\dots, d$ are supported on a proper subgroup $G\subset\z2$. Since
under our assumptions $a_r\neq b_r$ at least for one $r$, the subgroup is
two-dimensional, i.e. there exists a basis of $G$ consisting of two
vectors, say ${\bf u}$, $\bf v$. Let $(\alpha_r,\beta_r)$ be the
coordinates of $(a_r,b_r)$ related to this basis, i.e.
$(a_r,b_r)=\alpha_r {\bf u}+\beta_r{\bf v}$. Now define a new
two-dimensional random walk with steps

$$
\widehat Y_i=(\alpha_r,\beta_r) \quad\textrm{if}\quad
\widetilde Y_i=(a_r,b_r),\quad r=1,2,\dots, d,
$$
$$
\widehat Y_i=-(\alpha_r,\beta_r) \quad\textrm{if}\quad
\widetilde Y_i=-(a_r,b_r),\quad r=1,2,\dots, d.
$$
Then (cf. \cite{S76}, p. 65)
$$
\widehat Z_n=\widehat Y_1+\dots+\widehat Y_n,\quad
n=1,2,\dots
$$
is a symmetric recurrent aperiodic random walk in $\z2$. Obviously
$\widetilde Z_n={\bf 0}$ if and only if $\widehat Z_n={\bf 0}$. Hence we
have also $V^{S_{d-2}}(n)=\xi^{\widehat Z}(\bo, n)$ with obvious notation
for the local time of $\{\widehat Z\}$.

Now we prove the upper bound in Theorem 1.5, i.e.
\begg
\limsup_{n\to\infty}\frac{\xi_{S_{d-2}}^{(d)}(n)}{\log\log n}
\leq\lambda_d\has
\label{5up}
\endd

In the proof we follow the same lines as in the proof of Theorem
1.4. Choose $\lambda>\lambda_d$ and $\ep>0$ such that
$\delta=\lambda/\lambda_d-1-\ep>0$. Then using (\ref{geo}), Fact 5
and Fact 10  we can see as in (\ref{412}) that for $n$ big enough,
$$
\pe\left(\xi_{S_{d-2}}^{(d)}(n)>u\right)\leq (\log
n)^{1+\ep}P(\xi^{(d)}(\bo,n)\ge u) +\pe(V^{(d-2)}(n)>(\log
n)^{1+\ep})
$$
$$
\leq (\log n)^{1+\ep}\pe (\xi^{(d)}(\bo,\infty)>u)+
\pe(\xi^{\widehat Z}(\bo,n)>(\log n)^{1+\ep})
$$
$$
\leq ( \log n)^{1+\ep} \exp\left(-\frac{u} {\lambda_d}\right)+ \exp(-(\log
n)^{\ep/2}).
$$
 Hence choosing $n_j=\left[\exp\left(j^{2/\delta}\right)\right]$,
$$
\pe(\xi_{S_{d-2}}^{(d)}(n_{j+1})\ge\lambda\log\log n_j) \leq
C(\log n_j)^{-\delta}=\frac {C}{j^2}.
$$
Borel-Cantelli lemma and the usual monotonicity argument yields
(\ref{5up}), since $\varepsilon$ can be arbitrary small.

\medskip
Now we prove the lower bound
\begg
\liminf_{n\to\infty}\frac{\xi_{S_{d-2}}^{(d)}}{\log\log n}
\geq\lambda_d\has
\label{5low}
\endd

Here again we follow the proof of Theorem 1.4.
Let $\lambda<\lambda_d$ and $n_j=[\exp(j^{2/(1-\beta)})]$ with some
$\lambda/\lambda_d<\beta<1$ and define
\begg
\eta_0^j:=0,\qquad \eta_k^j:=\min\{i>\eta_{k-1}^j+[(\log n_j)^2]:\,
\bs_i\in S_{d-2}\},\qquad \by_k^j:=\bs_{\eta_k^j}.
\endd

Furthermore, let $\nu_j$ be the largest integer $N$ for which
\begg
\eta_N^j +[(\log n_j)^2] \leq n_j.
\endd

Since $\{\widehat Z\}$ is a two-dimensional random walk for which
Lemma 2.4 holds, choosing $a=[(\log n_j)^2]$, $f(n_j)=(\log
n_j)^\beta$ there, we get
$$
\pe(\nu_{j}\leq (\log n_j)^\beta)\leq \frac{c\log\log n_j}
{(\log n_j)^{1-\beta}}.
$$

Similarly to (\ref{gam100}), the probability
that a point will be visited at least $\lambda\log\log n_{j+1}$
times within $[(\log n_j)^{2}]$ steps is greater than
$$
\left(1-\gamma_d+O\left(\frac1{(\log n_j)^{1/2}}\right)\right)
^{[\lambda\log\log n_{j+1}]}=(\log
n_{j+1})^{-\lambda/\lambda_d}(1+o(1)).
$$

Define the events
$$
A_k^{j}=\{\xi^{(d)}(\bs_{\eta_k^{j}},\eta_k^{j}+[(\log n_j)^2])
-\xi^{(d)}(\bs_{\eta_k^{j}},\eta_k^{j})<\lambda\log \log
n_{j+1}\},
$$
then clearly  $A_k^{j}$ are independent in $k$ for fixed $j,$ and
have the same probability, hence we have for $j$ big enough
$$
\pe(A_k^{j})\leq 1-\frac{1+o(1)}{(\log
n_{j+1})^{\lambda/\lambda_d}}.
$$
Moreover,
$$
\pe(\xi_{S_{d-2}}^{(d)}(n_j)<\lambda\log\log n_{j+1})\leq
\left(1-\frac{1+o(1)}{(\log
n_{j+1})^{\lambda/\lambda_d}}\right)^{(\log
n_j)^\beta}+\frac{c\log\log n_j}{(\log n_j)^{1-\beta}}
$$
$$
\leq \exp\left(-(1/2)(\log n_j)^{\beta-\lambda/\lambda_d}\right)
+\frac{c\log\log n_j}{(\log n_j)^{1-\beta}}.
$$
This is summable in $j$, hence Borel--Cantelli lemma and monotonicity
imply (\ref{5low}). This together with (\ref{5up}) completes the proof
of Theorem 1.5. $\Box$

\bigskip\noindent
{\bf Proof of Theorem 1.6}. First we prove the lower bound, i.e.
\begg
2\lambda_d\leq \liminf_{n\to\infty}\frac{\xi_{B(r_n)}^{(d)}(n)}{\log
r_n}\has
\label{6lower}
\endd

For any $\delta>0$ and large enough $n$ we have $r_n^{2-\delta}\leq n$ and
by the law of the iterated logarithm (\ref{lil}) we conclude that for any
$\delta>0$, $\{\bs_i,\, 1\leq i\leq r_n^{2-\delta}\leq n\}$ are all in
$B(r_n)$ for large $n$ with probability one, hence by Theorem F we have
for any $\eta$ and large enough $n$,
$$
\xi_{B(r_n)}^{(d)}(n)\geq \xi^{(d)}(r_n^{2-\delta}) \geq
(1-\eta)\lambda_d\log r_n^{2-\delta}=
(1-\eta)(2-\delta)\lambda_d\log r_n \has
$$
Since $\eta$ and $\delta$ are arbitrary, (\ref{6lower}) follows.

To show the upper bound
\begg
 \limsup_{n\to\infty}\frac{\xi_{B(r_n)}^{(d)}(n)}{\log
r_n} \leq 2\lambda_d \has,
\label{6upper}
\endd
note that by (\ref{esc}) for any $\ep>0$ and large $n$ the random walk
does not hit points in the ball $B(r_n)$ after $r_n^{2+\ep}$ steps with
probability one. Hence by Theorem F for large $n$
$$
\xi_{B(r_n)}^{(d)}(n)\leq\xi^{(d)}(r_n^{2+\ep})\leq (2+2\ep)\ld \log r_n
\quad\has
$$
Since $\ep$ is arbitrary, (\ref{6upper}) follows.

This completes the proof of Theorem 1.6.
$\Box$

\renewcommand{\thesection}{\arabic{section}.}
\section{Further consequences and questions.}
\renewcommand{\thesection}{\arabic{section}}
\setcounter{equation}{0} \setcounter{theorem}{0}
\setcounter{lemma}{0}

{\bf 5.1 A special one dimensional walk.}

\bigskip

\nod Consider a simple symmetric walk in $\z2$
$$\bs_n^{(2)}=\sum_{k=1}^n \bx_k =(S_{n,1},S_{n,2})=(\sum_{k=1}^n
X_{k,1},\sum_{k=1}^n X_{k,2})$$

\nod Clearly the components $X_{k,1}$ and $ X_{k,2}$ are
dependent. However it is easy to check that the  pair

$$Y_{k,1}=X_{k,1}+X_{k,2},\qquad Y_{k,2}=X_{k,1}-X_{k,2}$$
are independent, with common distribution
 $$\pe(Y_{k,1}=\pm 1)=\pe(Y_{k,2}=\pm 1)=\frac{1}{2}.$$
Consequently,

$$V_n=\sum_{k=1}^n Y_{k,1}= S_{n,1}+S_{n,2} \quad {\rm and}\quad  Z_n
=\sum_{k=1}^n
Y_{k,2}=S_{n,1}-S_{n,2},\quad n=1,2 \dots$$ are independent simple
symmetric random walks.

Consider now the consecutive return times of $\{Z_n\}$ to zero, that
is to say let
$$\rho_0:=0, \quad\rho_k:=\min\{i>\rho_{k-1}:\, Z_i=0\},
\qquad k=1,2,\dots$$
Now clearly our two-dimensional walk $\bs_n^{(2)}$ is on the
line $x_1-x_2=0$ if an only if  $Z_n=0$, that is to say at the steps
$\rho_i,\,\, i=1,2, \dots $ Introducing the i.i.d. sequence
$U_i=V_{\rho_i}-V_{\rho_{i-1}}\,\,
i=1,2...$ we get a one-dimensional random walk

$$
R_n=V_{\rho_n}=\sum_{i=1}^n U_i, \qquad n=1,2,\dots
$$
For this walk we have from Spitzer \cite{S76}, p. 89 that
$$\pe(U_1=0)=1-\frac{2}{\pi},$$
$$\pe(U_1=2k)=\frac{2}{\pi}\frac{1}{4k^2-1}, \qquad k=\pm1,\pm2,
\dots $$

$R_n$ is in the domain of attraction of the Cauchy distribution,
so we will refer to it as Cauchy walk. Some properties of the
Cauchy walk was investigated in Taylor \cite{T99}. Here we want to point
out that our results have some implications for the local time of
the Cauchy walk. As
$$R_n=2k\Leftrightarrow\bs^{(2)}_{\rho_n}=(k,k),$$
we conclude that for the local time of $R_n$
\begg
\eta(2\ell,n):=\#\{i: 1\leq i \leq n, R_i=2\ell \}=\xi^{(2)}((\ell,
\ell),\rho_n). \label{ca1}
\endd
Thus from  (\ref{ca1}) we get that
$$\eta(n)=\max_y \eta(y,n)=\xi^{(2)}_L(\rho_n)$$
where $L=L(1,-1)$ is the line $x_1-x_2=0.$

Taking into account the well-known fact that $\log \rho_n\sim
2\log n$ (see e.g. \cite{R90}, p. 115) a simple application of our
Theorem 1.1 implies that for the maximal local time of the Cauchy
walk we have
$$ \frac{1}{2\pi}\leq\liminf_{\noo}\frac{\eta(n)}{(\log n)^2}\leq
\limsup_{\noo} \frac{\eta(n)}{(\log n)^2}\leq \frac{2} {\pi}\has
$$

\bigskip\noindent
{\bf 5.2 Open questions.}
\bigskip

1./ Our methods are not powerful enough to get exact constants in
Theorems 1.1, 1.2. We don't have any conjecture whether in
these theorems the $\liminf $ and $\limsup$ can be replaced by
limit and if so what would be the actual value of those limits. In
Theorem 1.3 even the exact order escapes us.

\bigskip

2./ In Theorems 1.1, 1.4 and 1.5 we have lines and subspaces going
through the origin. These results remain valid for lines and
subspaces having a fixed distance from the origin. However it would
be interesting to investigate the maximal local time on lines and
subspaces having a distance from the origin $d(n)\to \infty.$

\bigskip

3./ In our theorems we have balls centered at the origin. One
might be interested in the maximal local time in balls having a
center with distance $d(n)\to \infty$ from the  origin.

\bigskip

4./ We discussed subspaces and balls in the theorems but other
sets would be just as  interesting to be investigated. One
possibility would be to investigate angular domains, cones and
wedges. E.g. in Theorem 1.4 if we take the wedge in between two
planes, then we have to have a transition  from
$\frac{\lambda_d}{2}$ to $\lambda_d$ as the angle of the wedge
increases from $0\to 2\pi.$ Similarly we might ask that if in
Theorem 1.5 we replace the line with a cone (centered at the
origin) what kind of transition happen as the cone gets wider,
that is to say how does the order change from $\log\log n$  to
$\log n$ as the angle  of the cone changes form $0\to  \pi$ (from
the line to ${\cal Z}_3.$) Similar questions can be asked in two
dimension as well but those would be more interesting if the exact
constant was known in Theorem 1.1.

\bigskip

5./  One might ponder that how important is the actual shape of
the subset on which the maximum is taken. E.g. it would be
interesting to have results  on the maximal local time on sets
which are not specified in shape just their size are given and of
course one has to ensure that they are located close enough to the
origin that visits should occur.

\bigskip
6./ All our theorems are about simple symmetric walk. We might ask
how important is this restriction. Theorem C is valid  for a much
wider class of random walks.  So we might ask whether our theorems
remain valid for aperiodic   random walks under certain moment
conditions. In Theorem C all the moments has to exist, maybe
somewhat less precise results can be ensured even  if only the
second moment exists. Is it possible to say anything without
second moment?


\begin{thebibliography}{abcd}

\bibitem{A90} Auer, P.: The circle homogenously covered by
random walk on ${\cal Z}_2$. {\sl Statist. and Probab.
Lett.} \textbf{9} (1990), 403--407.

\bibitem{CH49} Chung, K.L. and Hunt, C.A.: (1949) On the
zeros of $\sum_1^n \pm 1.$ {\sl Annals of Math.} \textbf{50},
385--400.

\bibitem{CSF83} Cs\'aki, E. and F\"oldes, A.: How big are the increments
of the local time of a recurrent random walk? {\sl Z. Wahrsch. verw.
Gebiete} \textbf{65} (1983), 307--322.

\bibitem{CSRR98} Cs\'aki, E., R\'ev\'esz, P. and Rosen, J.:
Functional laws of the iterated logarithm for local times of
recurrent random walks  on $\z2$. {\sl  Ann. Inst. H. Poincar\'e
Probab. Statist.} \textbf{34} (1998), 545--563.

\bibitem{DPRZ01}  Dembo, A., Peres, Y., Rosen, J. and
Zeitouni, O.: Thick points for planar Brownian motion and the
Erd\H{o}s-Taylor conjecture on random walk. {\sl Acta Mathematica
} \textbf{186} (2001), 1--30.

\bibitem{DE50} Dvoretzky, A. and Erd\H os, P.: Some problems on random
walk in space. {\sl Proc. Second Berkeley Symposium} (1950), 353--368.

\bibitem{ET60} Erd\H os, P. and Taylor, S.J.: Some problems
concerning the structure of random walk paths. {\sl Acta Math. Acad.
Sci. Hung.} \textbf{11} (1960), 137--162.


\bibitem{MR94} Marcus, M. and Rosen, J.: Laws of
the iterated logarithm for the local times of recurrent random walks
on $\z2$ and of L\'evy processes and random walks in
the domain of attraction of Cauchy random variables. {\sl Ann.
Inst. H. Poincar\'e Prob. Stat.} \textbf{30} (1994), 467--499.

\bibitem{R90} R\'ev\'esz, P.: {\sl Random Walk in
Random and Non-Random Environments}. World Scientific, Singapore,
1990.

\bibitem{S76} Spitzer, F.: {\sl Principles of Random
Walk}, 2nd. ed. Van Nostrand, Princeton, 1976.

\bibitem{T99} Taylor, H.M.: The fundamental matrix for a certain
random walk. {\sl J. Appl. Probab.} \textbf{36} (1999), 320--333.
\end{thebibliography}
\end{document}